%Paper: alg-geom/9408003
%From: Oliver Kuechle <kuechle@math.ucla.edu>
%Date: Tue, 9 Aug 1994 12:54:08 -0700 (PDT)

%%%%%%%%%%%%%%%%%%%%%%%%%%%%%%%%%%%%%%%%%%%%%%%%%%%%%%%%%%%%
% LOCAL POSITIVITY OF AMPLE LINE BUNDLES
%
% by Lawrence Ein, Oliver K\"uchle and Robert Lazarsfeld
%
% AMS-TeX version 2.1, 23 pages
%%%%%%%%%%%%%%%%%%%%%%%%%%%%%%%%%%%%%%%%%%%%%%%%%%%%%%%%%%%%
\documentstyle{amsppt}
\magnification \magstep1
\parskip 10pt
\parindent .3in
\pagewidth{5.2in}
\pageheight{7.2 in}
\NoRunningHeads

\def \bl{\vskip 10pt}

\def \ni{\noindent}
\def \P{\bold{P}}
\def \Q{\bold{Q}}
\def \Z{\bold{Z}}
\def \R{\bold{R}}
\def \C{\bold{C}}
\def \O{\Cal{O}}
\def \D{\Cal{D}}
\def \lra{\longrightarrow}
\def \d {\delta}
\def \eps{\epsilon}
\def \mult{ \text{\rm mult}}
\def \pr {\prime}
\def \dim {\text{  dim} \phantom{.}}
\def \F {\Cal F}
\def \I {\Cal I}
\def \B {\Cal B}
\def \V {\Cal V}

\def\a  {\alpha}
\def\g  {\gamma}
\def\cv  {\cdot}
\vskip .8 in
\centerline{\bf LOCAL POSITIVITY OF AMPLE LINE BUNDLES}
\bl
\bl
\centerline{\smc Lawrence Ein \footnote{Partially supported by NSF Grant
DMS 93-02512}}
\centerline{\smc Oliver K\"uchle \footnote{Supported  by
Deutsche Forschungsgemeinschaft}}
\centerline{\smc Robert Lazarsfeld \footnote{Partially supported by
NSF Grant DMS 94-00815}}

\bl

\ni {\bf Introduction}

The purpose of this paper is to establish a  lower bound on
the Seshadri constants measuring the local positivity of an ample line bundle
at a general point of a complex projective variety of arbitrary dimension.

Let $X$ be an irreducible complex projective variety, and let $L$ be a nef line
bundle on $X$. Demailly \cite{De2} has introduced
a very interesting invariant
which in effect measures how positive $L$ is locally near
a given smooth point $x \in  X$. This {\it Seshadri constant} $\eps(L,x) \in
\R$  may be defined as follows. Consider the blowing up
$$f : Y = \text{Bl}_x(X) \lra X$$
 of $X$ at $x$, and denote by  $E = f^{-1}(x) \subset Y$ the exceptional
divisor.  Then $f^* L$ is a nef line bundle on $Y$, and we put
$$ \eps(L, x) =
\sup \left \{ \eps \ge 0 \ \mid  \ f^*L - \eps \cdot E \ \text{is nef}
 \ \right \}.$$
Here $\ f^*L - \eps E$ is considered as an $\R$-divisor on $Y$, and to say that
it is nef means simply that $f^*L \cdot C^\prime \ge \eps E \cdot C^\prime$ for
every irreducible curve $C^ \prime \subset Y$.  For example, if $L$ is
very ample, then $\eps(L, x) \ge 1$ for
every smooth point $x \in X$. Seshadri's criterion
(cf. \cite{Ha}, Chapter 1)
states that $L$ is
ample if and only if there is a positive number $\eps > 0$ such that $\eps(L,
x)
\ge \eps$ for every $x \in X$. We refer to \S1 below, as well as
\cite{De2}, \S 6, for alternative
characterizations and further properties of  Seshadri
constants.

It was shown by an elementary argument in \cite{EL} that
if $S$ is a smooth
projective surface, and $L$ is an ample line bundle on $S$,
then $\eps(L, x) \ge 1$ for all except perhaps countably many $x \in S$.
This suggested the somewhat
surprising possibility that there could be a similar lower bound on the
local positivity of an ample line bundle at a general point of an irreducible
projective variety of any dimension. Our main result shows that this is
indeed the case:
\proclaim{\bf Theorem 1}  Let $L$ be a nef and big line bundle on an
irreducible projective variety $X$ of dimension $n$.  Then  $$\eps(L,x)
\ge
\frac{1}{n}$$ for all $x \in X$ outside a countable union of proper closed
subvarieties of $X$.  Moreover if $L$ is ample, then given any $\delta >
0$ the locus
$$\left \{ x \in X \  \bigg   | \ \eps(L,x) > \frac{1}{n+\delta} \ \right
\}$$
contains a Zariski-open dense set. \endproclaim
\ni More generally, we prove  that if there exists a countable union $\B
\subset X$ of proper closed subvarieties,  plus a real number $\alpha > 0$
such that  for $1 \le r \le n$:
$$\int_Y c_1(L)^r  \ge (r \cdot \alpha)^r \ \  \ \ \forall \
\ \text{$r$-dimensional } Y \subset X \ \text{ with } Y \not \subset \B,$$
 then $\eps(L,x) \ge \alpha$ for all sufficiently
general $x \in X$. Examples constructed by Miranda show that
given any $ b > 0$, there exist $X , L  \text{ and }x$ such that $0 <
\eps(L, x) < b$. In other words, there cannot be a bound (independent
of $X$ and $L$) that holds at every point. On the other hand, it is
unlikely that the particular constant appearing in Theorem 1 is
optimal.   In fact, it is natural to conjecture that in the setting of
the Theorem one should have $\eps(L, x) \ge 1$ for a very general point
$x \in X$.

Recent interest in Seshadri constants stems in part from the fact that
they govern an elementary method for producing sections of adjoint
bundles. Our bounds then imply the following, which complements the
non-vanishing theorems of Koll\'ar (\cite{Ko1} \S3):
\proclaim{\bf Corollary 2} Let $L$ be a nef line bundle on a smooth
projective variety $X$ of dimension $n \ge 2$, and given an integer $s \ge
0$ suppose that $$\int_Y c_1(L)^r \ge (r(n+s))^r$$ for every
$r$-dimensional subvariety $Y \subseteq X$ not contained in some fixed
countable union $\B \subset X$ of proper subvarieties. Then the adjoint
series $|\O_X(K_X + L )|$ generates $s$-jets at a general
point $x \in X$, i.e. the evaluation map
$$H^0(X,\O_X(K_X + L ) ) \lra H^0(X, \O_X(K_X + L) \otimes \O_X/
\I_x^{s+1})$$ is surjective. In particular,
$$h^0(X, \O_X(K_X + L)) \ge \binom{n+s}{n}.$$ \endproclaim
\ni It follows for example that if $A$ is ample, then  $\O_X(K_X +
(ns + n^2 )A)$ generates $s$-jets at almost all points $x \in X$.
We remark that contrary to what one might expect from extrapolating the
well-known conjectures of Fujita \cite{Fu} on global generation and
very ampleness, there cannot exist a linear function $f(s)$
(depending on $n$,  but independent of $X$ and $A$) such that
 $\O_X(K_X + f(s)A)$ generates $s$-jets for $s \gg 0$ at  {\it every}
point of $X$ (Remark 1.7).

Similarly:
\proclaim{Corollary 3}  Suppose that $L$ is a  nef and big line bundle on
a smooth projective variety $X$ of dimension $n \ge 2$.  Then for all
$m \ge 2n^2$, the linear series $|K_X + mL|$ is birationally very
ample, i.e. the corresponding rational map $$\phi_{|K_X + mL|} : X
\dashrightarrow \P$$ maps $X$ birationally onto its image. \endproclaim
\ni For example, suppose that $X$ is a smooth minimal variety of
general type, i.e. $K_X$ is nef and big. Then the pluricanonical
rational maps
$$\phi_{|mK|} : X \dashrightarrow \P$$
are birational onto their images
for $m > 2n^2 $. This extends (with somewhat weaker numbers)  results
of Ando \cite{An}  in the cases $n \le 5$.  More generally, if $X$ is
a general type minimal $n$-fold of global index $r$, then $|mrK_X|$ is
again birationally very ample when $m > 2n^2$ (Corollary 4.6). As above,
one also  has an analogue of Corollary 3 for the linear series $|K_X + L|$
involving intersection numbers of $L$ with subvarieties of $X$.

The proof of Theorem 1 draws  inspiration from two sources:  first, the
arguments used in \cite{Nad}, \cite{Ca} and
\cite{KoMM} to prove boundedness of Fano manifolds of Picard number
one; and secondly,  some of the geometric ideas occuring in
\cite{Be}, \cite{Nak} and especialy \cite{Fa}.
Roughly speaking, if Theorem 1 fails then given a general point $x \in X$
there exists a curve $C_x
\subset X$ through $x$ such that
$$\frac{\text{mult}_x(C_x)}{  (L \cdot C_x)} > n.$$ We start by fixing a
divisor $E_x \in |kL| \text{ for } k \gg 0$  with suitably large
multiplicity at $x$. If we could  arrange that $C_x \not
\subset E_x$,  then one  arrives right away at a contradiction by
estimating $E_x \cdot C_x$ in terms of multiplicities at $x$.
Unfortunately it doesn't seem to be immediate that one can do so. Instead, we
use a gap construction to show that for an appropriate choice of $y = y(x)$,
we can at least control the difference of  the multiplicities of $E_x$ at $y$
and at a general point of $C_y$. The principal new ingredient is then an
argument showing that we can {\it rechoose} the divisors $E_x$  in such a
way as to ensure that
$C_y \not \subset E_x$ while keeping the multiplicity mult$_{y}(E_x)$
of $E_x$ at $y$ fairly large, and then we are
done. Stated somewhat informally, the main lemma here is the following:

\pagewidth{4.3in}
\parindent .6in
\item
\ni Suppose that $\{ Z_t \subseteq  V_t \}_{t\in T} $ is a family of
subvarieties of a  smooth variety $X$, parametrized by a smooth affine
variety $T$, and assume that $\cup_{t \in T} V_t$ is dense in $X$. Suppose
also given a family   $\{ E_t\}_{t\in T} \in |L| $ of divisors in a fixed
linear series on $X$, with
$$a = \mult_{Z_t}(E_t)  \ \text{ and } \ b = \mult_{V_t}(E_t)$$ for
general $t \in T$. Then one can find another family of divisors
$\{ E^\pr_t \}_{t \in T} \in |L| $  such that
$$V_t \nsubseteq E^\pr_t  \ \text{ and }  \ \mult_{Z_t}(E_t^\prime)
\ge a - b$$ for general $t\in T$.

\pagewidth{5.2in}
\parindent .3in
\ni We refer to Proposition 2.3 for the precise statement and proof.
Denoting by \{$s_t\} \in \Gamma(X, L)$  the family of sections
defining $E_t$, the idea is  to construct $E_t^\prime$ as the
divisor of the section $D s_t \in  \Gamma(X, L)$, where $D$ is a general
differential operator of order $ b$ in $t$. For
divisors on projective space (and other compactifications of group
varieties), the process of differentiating in order to arrive at a proper
intersection plays an important role e.g. in \cite{Be} and
\cite{Fa}. Our observation here is that the same idea works in a
deformation-theoretic context, and we hope that the lemma may find other
applications in the future.

Concerning the organization of the paper, we start in \S 1 with a quick
review of general facts about Seshadri constants. In \S 2 we discuss
multiplicities in a family of divisors, and show in particular that by
differentiating in parameter directions one can lower the multiplicity of
such a family along a covering family of subvarieties. The proof of the
main result occupies \S 3, and finally we give some elementary
applications in \S 4.

We have profitted from discussions with F. Campana, V.
Ma{\c s}ek, A. Nadel and G. Xu. Nadel in particular stressed  some
years ago the relevance of techniques from diophantine approximation
and transcendence theory to arguments of this type. We are especially
endebted to M. Nakamaye for helping us to understand some of these
arithmetically motivated ideas. In particular, the crucial Proposition
2.3 was inspired by a proof Nakamaye showed us of Dyson's lemma
concerning singularities of curves in $\P^1 \times \P^1$.
\bl

\ni{\bf \S 0. Notation and  Conventions.}

(0.1). We work throughout over the complex numbers $\C$.

(0.2). We will say that a property holds at a {\it general} point of a
variety $X$ if it holds for a non-empty Zariski-open subset of $X$. It
holds at a {\it very general} point if it is satisfied off the union of
countably many proper closed subvarieties of $X$.

(0.3).  If $X$ is a projective variety of dimension $n$, and $L$ is a line
bundle on $X$, we denote by $L^n \in \Z$ the top self-intersection
number of $L$. Given a subvariety $Z \subset X$ of dimension $r$, $L^r
\cdot Z$ indicates the degree $\int_Z c_1(L)^r \in  \Z$.   Recall that a line
bundle $L$ is {\it numerically effective}
or {\it nef} if $L \cdot C \ge 0$ for all effective curves $C \subseteq X$.
Kleiman's criterion (\cite{Ha}, Chapter1) implies that a line
bundle is nef if and only if it lies in the closure of the ample cone in the
N\'eron - Severi vector space $NS(X)_\R$. Recall also that a nef line
bundle $L$ is big if and only if $L^n > 0$. Similar definitions and
remarks hold for (numerical equivalence classes of)
$\Q$-Cartier $\Q$-divisors on $X$.

(0.4). For  varieties $X$ and $T$, $pr_1 :X \times T \lra X, \ \ pr_2 : X
\times T \lra T$ denote the projections. If $Z \lra T$ is a mapping,
$Z_t$ denotes the fibre of $Z$ over $t \in T$. Given a Zariski-closed subset
(or subscheme) $Z \subset X \times T$, we consider the fibre $Z_t$ of $pr_2$
as a subset (or subscheme) of $X$. Similarly, $Z_x \subset T$ is the
fibre of $Z$ over $x \in X$. If $V \subset X$ is a subvariety, $\I_V
\subset \O_X$ denotes its ideal sheaf.

\vskip 14 pt

\ni{\bf \S 1. Seshadri Constants.}
\vskip 3 pt

In this section, we recall briefly some of the basic facts about
Seshadri constants. We start with a

\subhead
Definition 1.1
\endsubhead
(\cite{De2}) Let $X$ be a  projective
variety, $x\in X$ a point and $L$ a nef line bundle on $X$.
Then the {\sl Seshadri constant of $L$ at $x$} is the real number
$$\eps (L,x):= \inf_{C\owns x} \Biggl\{{L\cv C\over
\text{mult}_x(C)}\Biggr\},$$
where the infimum is taken over all reduced and irreducible curves
$C\subset X$ passing through $x$. (We remark that it is enough here that
$L$ be a $\Q$-Cartier $\Q$-divisor.)

It is elementary (and standard) that (1.1) is equivalent to the
alternative definition given in the Introduction:
\proclaim{Lemma 1.2}
Let $X$ be a projective variety, $L$ a nef line bundle on $X$, and $x
\in X$ a smooth point. Let
$$f : Y = \text{Bl}_x(X) \lra X$$
be the blowing up of $X$ at $x$, and denote by  $E = f^{-1}(x) \subset Y$
the exceptional divisor.  Then
$$ \eps(L, x) = \sup \left \{ \eps \ge 0 \ \mid
\ f^*L - \eps \cdot E \ \text{is nef} \ \right \}. \ \qed $$\endproclaim
\ni Note that the supremum in the definition is actually a maximum.

There are interesting characterizations of Seshadri constants involving the
generation of jets. Recall that  given  a line bundle
$B$ on a smooth variety $X$, and an integer $s\ge 0$, we say that the linear
series $|B|$ {\sl generates s-jets at $x\in X$} if the  evaluation
map
$$H^0(X,\O_X(B) ) \lra H^0(X, \O_X(B) \otimes \O_X/ \I_x^{s+1})$$
is surjective, where $\I_x$ denotes the ideal sheaf of $x$. The following
Proposition --- which is a variant of \cite{De2}, Theorem 6.4 ---
shows in effect that computing the Seshadri constant $\eps(L,x)$ is
equivalent to finding a linear function $f(s)$ such that the adjoint
series  $|K_X + f(s)L|$ generates $s$-jets at $x$ for all $s \gg 0$.
\proclaim{ Proposition 1.3} Let $L$ be a nef and big line bundle on a
smooth projective variety $X$ of dimension $n$.

\parskip 5pt
\ni {\rm (1.3.1).}
If $${\displaystyle r > {s\over \eps(L,x)} + {n\over
\eps(L,x)}},$$  then $|K_X+rL|$ generates $s$-jets at $x\in X$. The same
statement holds if
 ${\displaystyle r = {s +n\over \eps(L,x)}}$ and $L^n>\eps(L,x)^n$.

\ni {\rm (1.3.2).}
Conversely, suppose there is a real number $\eps >0$ plus a constant $c \in
\R$ such that $|K_X+rL|$ generates $s$-jets at $x\in X$ for all $s\gg 0$
whenever
$${\displaystyle  r > {s\over \eps}+c}.$$  Then $
\eps(L,x)\ge \eps.$
\endproclaim
\parskip 10pt

\demo{\bf Proof}
(1). (\cite{De2}, Prop. 6.8.)  This is a standard application of
Kawamata-Viehweg vanishing for nef and big line bundles. In brief, let
$f:Y = \text{Bl}_x(X) \lra X$ be the blowing-up of
$X$ at $x$, with exceptional divisor $E \subset Y$. It suffices to
show that
$$H^1(X,\O_X(K_X+rL)\otimes \I_x^{s+1})=H^1(Y,\O_X(K_Y+rf^*L-(s+n)E)) = 0.
\tag *$$ Setting $\eps=\eps(L,x)$, one has the numerical equivalence
$$rf^*L -(s+n)E\equiv {s+n\over \eps}(f^*L-\eps E)+(r-{s+n\over
\eps})f^*L.$$ Hence $r f^*L - (s+n)E$ is big and nef, and (*) follows from
Vanishing.
\par
(2). Let $C\owns x$ be a reduced and irreducible curve with mult$_x(C)=m$.
Fix $s\gg 0$ and let $r$ be the least integer $>(s/\eps) +c$.
The geometrical interpretation of the fact that $|K_X+rL|$ generates
$s$-jets at $x$ is that we can find a divisor $D_x\in |K_X+rL|$, with
mult$_x(D_x)=s$,  having an arbitrarily prescribed tangent cone at $x$.
In particular, we can choose $D_x$ such that the tangent cones to $D_x$ and
$C$ at $x$ meet properly, and  since $C$ is irreducible it follows that
$D_x$ and $C$ themselves meet properly.  Then
$$C\cv (K_X+rL)\ge \text{mult}_x(C)\cdot\text{mult}_x(D_x)=m\cdot s,$$
and hence
$${C\cv L\over m}\ge {s\over r}-{C\cv K_X\over rm}.$$
Since $r\le (s/\eps) +c +1$ the claim  follows by
letting $s\to\infty$.
\qed \enddemo

Our main result (Theorem 3.1) will give a lower bound on  the Seshadri
constant of a nef and big line bundle at a very general point, i.e. a bound
which holds off the union of countably many proper subvarieties. However
the following Lemma shows that at least for ample line bundles, one then
obtains a statement valid on a Zariski-open set:
\proclaim{Lemma 1.4}  Let $L$ be an ample line bundle on an irreducible
projective
variety $X$. Suppose that there is a positive rational number $B > 0$ and a
smooth point $y \in X$ for which one knows that $\eps(L, y) > B$. Then the
locus
$$\left \{ x \in X \ \bigg | \ \eps(L,x) > B \right \}$$
contains a Zariski-open dense set. \endproclaim
\demo{\bf Proof} Given a smooth point $x \in X$, let
$$f_x : Y_x = \text{Bl}_x(X) \lra X$$
be the blowing-up of $X$ at $x$, with exceptional divisor $E_x \subset
Y_x$. Consider the $\Q$-divisor
$$M_x =_{\text{def}} f_x^*L - B \cdot E_x$$
on $Y_x$. It follows from the hypothesis and characterization (1.2) of
Seshadri constants that  $f_y^*L - \eps(L,y) \cdot E_y$ is nef, and
hence $M_y$ is ample. Since ampleness is an open condition in a flat family
of line bundles, there exists a non-empty Zariski-open subset $U \subset X$
of smooth points of $X$ such that $M_x$ is ample whenever $x \in U$. But by
(1.2) again, the ampleness of $M_x$ implies that $\eps(L,x) > B$.   \qed
\enddemo

Finally, for the convenience of the reader we recall from
\cite{EL}, \S3, the examples of Miranda showing that the Seshadri
constants of an ample line bundle can take on arbitrarily small positive
values.
\proclaim{Proposition 1.5} {\rm  ({\bf Miranda})} Given any positive number
$b > 0$, there exist a projective variety $X$ and an ample line bundle $L$ on
$X$ such that
$$0 < \eps(L, x) < b$$
for all $x$ in a codimension two subset $V \subset X$. \endproclaim
\demo{\bf Sketch of Proof} We first construct  a line bundle $N$ on a
surface $S$ having the required property. To this end, start with a reduced
and irreducible plane curve $C \subset \P^2$ of degree $d \gg 0$ with a
point $y \in C$ of multiplicity $m > \frac{1}{b}$. Fix  a second integral
curve $C^\prime \subset \P^2$ of degree $d$ meeting $C$ transversely.
Provided that $d$ is sufficiently large, we can assume by taking
$C^\prime$ generally enough that all the curves in the pencil spanned by
$C$ and $C^\prime$ are reduced and irreducible. Blow up the base-points of this
pencil to obtain a
surface $S$, admitting a map $f : S \lra \P^1$ with irreducible fibres,
among them $C \subset S$. Observe that any of the exceptional divisors
over $\P^2$ gives rise to a section $\Gamma \subset S$ of $f$ which meets
$C$ transversely at one point. Fix an integer $a \ge 2$, and put $ N = aC
+ \Gamma$. It follows from the Nakai criterion that $N$ is ample. But $N
\cdot C = 1$ whereas $\mult_y(C) = m > \frac{1}{b}$, so $\eps(N, y) < b$.
As Viehweg pointed out, this gives rise automatically to  higher
dimensional examples. In fact, take for instance $X = S \times
\P^{n-2}$ and $L = pr_1^*(N) \otimes pr_2^*(\O_{\P}(1))$.  By considering
the evident curve in $S \times \{ z \}$, one sees that
$$\eps(L, (y,z)) \le \eps(S,y) < b \ \ \text{for all } \ z \in \P^{n-2}.$$
Thus it suffices to take $V = \{ y \} \times \P^{n-2}$. \qed
\enddemo

\ni {\bf Remark 1.6.}  We do not know whether Seshadri constants can become
arbitrarily small on a codimension one subset of $X$. It is shown in
\cite{EL} that this cannot happen when $X$ is a surface.

\ni{\bf Remark 1.7.}  A well-known conjecture of Fujita \cite{Fu}
asserts that if $L$ is an ample line bundle on a smooth projective variety
$X$ of dimension $n$, then $\O_X(K_X + (n+1)L)$ is free and $\O_X(K_X +
(n+2)L)$ is very ample. Extrapolating, one might be tempted to hope that
for all $s \ge 0$:
$$|K_X + (n+s+1)L| \ \text{ separates $s$-jets at every point } x \in
X. \tag * $$
However  (1.3.2) and (1.5) show that (*) is not true in general. In fact,
there cannot exist a linear function $f(s)$ (depending on $n$ but
independent of
$X$ and $L$) such that $|K_X + f(s)L|$ generates
$s$-jets at all $x \in X$. However when $X$ is a surface, it follows from
\cite{BS} or \cite{De1} that there exists a
{\it quadratic} function $f(s)$ such that $|K_X + f(s)L|$ separates $s$-jets
at every $x \in X$. (Cf. \cite{La}, \S7.)

\ni{\bf  Remark 1.8.}  Let $X$ be a projective variety of dimension $n$,
$L$ a nef line bundle on $X$, and $x \in X$ a smooth point. Then
$$\eps(L, x) \le \root n \of{ (L^n) } \ \text{ for every } x \in X. $$ In
fact, if $f : Y = \text{Bl}_x(X) \lra X$ is the blowing up of
$x$, with exceptional divisor $E$, then $(f^*L - \eps(L,x) \cdot E)^n \ge
0$. As an interesting example, suppose that $X$ is a simple abelian
variety, and  $L$ is a principal polarization on $X$. Then Nakamaye
has shown  that $\eps(L, x) \ge 1$ for all $x$. Therefore one has the
inequality:
$$1 \le \eps(L, x) \le \root n \of{n!} \approx \frac{n}{e}.  \tag 1.8.1
$$
Since $X$ is homogeneous, $\eps(L,x)$ is independent of $x \in X$, so
there is a real number $\eps(L)$ satisfying (1.8.1) canonically attached
to a principally polarized abelian variety $(X, L)$.  However it is not
obvious to us what the value of this invariant is, even for Jacobians or
very general p.p.a.v.'s. Note that if $C \subset X$ is any curve, then $L
\cdot C \ge n $. Hence (1.8.1) implies that the curves computing
$\eps(L,x)$ in (1.1) cannot be smooth.

\vskip 8 pt
\ni{\bf \S2.  Multiplicity Lemmas.}

This section is devoted to some preliminary results concerning
multiplicity loci in a family of divisors. Proposition 2.3 --- which
allows one to reduce the multiplicity of a family of divisors along
a covering family of subvarieties --- is the crucial ingredient in the
proof of our main Theorem. It is in this section that we make essential use
of the fact that we are working in characteristic zero.

We start with some notation. If $M$ is a smooth variety, and $E$ is an
effective divisor on $M$, then the function $x \mapsto \mult_x(E)$ is
Zariski upper-semicontinuous on $M$. Given an irreducible subvariety $Z
\subset M$, by $\mult_Z(E)$ we mean the value of $\mult_x(E)$ at a
general point $x \in Z$. We refer to (0.4) for notation and conventions
concerning projections from products, and fibres of morphisms.

The first lemma allows one to make fibrewise calculations of
multiplicities. It is certainly a well-known fact, but we include a
 proof for the convenience of the reader.
\proclaim {Lemma 2.1}  Let $X$ and $T$ be smooth irreducible
varieties, and suppose that $Z \subset X \times T$ is an irreducible
subvariety which dominates $T$ (under projection to the second factor).
Let $E \subset X \times T$ be any effective divisor. Then for a general
point $t \in T$, and any irreducible component $W_t \subseteq Z_t$ of the
fibre $Z_t$, we have:
$$
\mult_{W_t}(E_t) = \mult_Z(E).
$$
\endproclaim
\demo{\bf Proof} Consider more generally a
 mapping $f : M \lra T$ of smooth varieties, and suppose that $V \subset M$
is a smooth subvariety dominating $T$.  Assume given an effective divisor
$E \subset M$ with
$\mult_V(E) = a$. We will show that for a general point $t \in
T$, and any irreducible component $W_t \subseteq V_t$:
$$\mult_{W_t}(E_t) = a. \tag *$$
The Lemma then follows by taking $M$ to be an open subset
of $X \times T$ on which $Z$ is smooth, and setting $V = Z \cap M$.

To prove (*), note first that for	$b \le a$ the section
$s \in \Gamma(M, \O_M(E))$ defining $E$ lies in the subspace
$$\Gamma(M, \O_M(E) \otimes \I_V^b) \subseteq \Gamma(M, \O_M(E)).$$
 Hence $s$ determines a section
 $$\delta_b( s) \in \Gamma(V, \I_V^b/\I_V^{b+1}(E))
= \Gamma(V,Sym^b(N_{V/M}^*)(E))$$
of a twist of the $b^{\text{th}}$ symmetric power of the
conormal bundle to $V$ in $M$. (One thinks of
$\delta_b(s)$ as giving the $b^{\text{th}}$ order terms in the Taylor
expansion of $s$ in directions normal to $V$.) One checks e.g. by a
calculation in local coordinates that $\delta_b(s) = 0$ for $b < a$ whereas
$\delta_a(s) \ne 0$.  Now fix a point $t \in T$ lying in the open subset of
$T$ over which the mappings $M \lra T$ and  $V \lra T$ are smooth, and
let $s_t = s | M_t \in \Gamma(M_t, \O_{M_t}(E_t))
$ be the restriction of $s$ to $M_t$, so that $s_t$ is the section defining
$E_t$. Then
 $$\delta_b(s) | V_t  = \delta_b(s_t)  \in \Gamma(V_t,
Sym^b(N^*_{V_t/M_t})(E_t)).$$
But since $V \lra T$ is dominating, a non-zero section of a locally free
sheaf on $V$ restricts to a non-zero section on each irreducible component
$W_t \subseteq V_t$ of a general fibre. Hence $\delta_b(s_t) = 0$ for $b <
a$ and $\delta_a(s_t) \ne 0 \in \Gamma(W_t,\O_{W_t}(E_t))$ for
general $t \in T$. But as we have just seen, this implies that
$\mult_{W_t}(E_t) = a$, as claimed. \qed \enddemo

\ni{\bf Remark 2.2.} Some readers may prefer to see the argument phrased
in a  more concrete manner. In the situation of (2.1) it is enough
to show that for  sufficiently general  $t \in T$, and for any component
$W_t \subset Z_t$, there exists at least one point $x \in W_t$ such that
$\mult_x(E_t) = \mult_{(x,t)}(E)$. Now since $Z$ dominates $T$, given
general points $t  \in T$ and $x \in W_t \subset Z_t$ we can find a local
analytic section of the projection $Z \lra T$, say $\sigma : U \lra Z$,
defined in a (classical) neighborhood $U$ of $t$ in $T$, whose image
passes through the point $(x,t)$. Replacing $T$ by $U$, and working
analytically, we can assume given a holomorphic mapping
$p : T \lra X$, and we are reduced to proving that
$$\mult_{p(t)}(E_t) = \mult_{(p(t),t)}(E) \tag * $$
for general $t \in T$. But this follows easily from an explicit calculation
in local holomorphic coordinates. [Choose coordinates $x$ and $t$ on $X$
and $T$, and suppose $p$ is given by $p = p(t)$. Defining $y = x - p(t)$,
expand a local equation for $E$ as a Taylor series in $y$ and $t$.]

\vskip 5pt

We now come to the main result of this section.
\proclaim{Proposition 2.3}  Let $X$ and $T$ be smooth irreducible
varieties, with $T$ affine, and suppose that
$$Z \subseteq V \subseteq X \times T$$
are irreducible subvarieties such that $V$ dominates $X$. Let $L$ be a
line bundle on $X$, and suppose given on $X \times T$ a divisor
$$E \in | pr_1^*(L)|.$$
Write
$$\ell = \mult_Z(E), \ \ k = \mult_V(E).$$
Then there exists a divisor $$E^\pr \in |pr_1^*(L)|$$
on $X \times T$ having the property that
$$ \mult_Z(E^\pr) \ge  \ell  - k, \ \ \text{and } \ V \nsubseteq \text{
Supp}(E^\pr).$$
\endproclaim
\ni Let $\sigma \in \Gamma(X \times T, pr_1^*(L))$ be the
section defining
$E$. In a word, the plan is to obtain $E^\pr$ as the divisor of a
section
$$\sigma^\pr = D \sigma \in \Gamma(X \times T, pr_1^*(L)),$$
where $D$ is a general differential operator of order $\le k$ on $T$.
 So we begin with some remarks about
differentiating sections of the line bundle $pr_1^*(L)$ in parameter
directions.

	Let $\D_T^k$ be the (locally free) sheaf of differential operators of
order $\le k$ on $T$. Then sections of $\D_T^k$ act naturally on the
space $\Gamma(X \times T, pr_1^*(L))$ of  sections of
$pr_1^* (L)$. Naively this comes about as follows. Choose local
coordinates $x$ and $t$ on
$X$ and $T$, and let $g_{\alpha, \beta}(x) $ be the transition
functions of $L$ with respect to a suitable open covering of
$X$. Then sections of $pr_1^*(L) $ are given by collections of
functions $\sigma = \{ s_{\alpha}(x, t) \}$ such that
$ s_{\alpha}(x,t) = g_{\alpha, \beta}(x) s_{\beta}(x,t).$ If
$D$ is a differential operator in the $t$-variables, then
$$Ds_{\alpha}(x,t) = g_{\alpha, \beta}(x) D s_{\beta}(x,t).$$
Therefore the $\{ D s_{\alpha}(x,t) \}$ patch together to define a section $D
\sigma \in \Gamma( X \times T , pr_1^*(L))$.

	To say the same thing in a more invariant fashion, let
$\D^k_{X\times T}(pr_1^*(L))$ denote the sheaf of differential
operators of order $\le k$ on $pr_1^*(L)$, i.e.
$$
\D^k_{X\times T}(pr_1^*(L)) = P^k_{X \times T}(pr_1^*(L)) ^*
\otimes pr_1^*(L),
$$
where $ P^k_{X \times T}(pr_1^*(L)) $ is the
sheaf of principal parts associated to $pr_1^*(L)$. Observe that there
is a canonical inclusion of vector bundles
$$pr_2^* (\D^k_T) \hookrightarrow \D^k_{X \times T}(pr_1^*(L)). \tag
*
$$
In fact, it follows from the construction of bundles of principal
parts plus the projection formula that one has an isomorphism:
$$P^k_{X \times T   /  X}(pr_1^*(L)) =
pr_2^*(P^k_T(\O_T)) \otimes pr_1^*(L),$$
and then (*) is deduced from the surjection
$$P^k_{X \times T}(pr_1^*(L)) \lra P^k_{  X \times T
/  X}(pr_1^*(L)).$$
On the other hand, a section $\sigma \in \Gamma(X \times
T, pr_1^*(L))$ gives rise to a vector bundle map
$$ \D^k_{X\times T}(pr_1^*(L)) \lra pr_1^*(L),$$
and hence by
composition a homomorphism
$$ j_{\sigma} : pr_2^* (\D^k_T) \lra pr_1^*(L).$$
 Given
$D \in \Gamma(T, \D^k_T )$, $$D \sigma \in \Gamma(X \times T, pr_1^*(L))$$ is
just the image of
$pr_2^*(D) \in \Gamma(X \times T, pr_2^* (\D^k_T))$ under the map
on sections determined   by $j_{\sigma}$.

\demo{\bf Proof of Proposition 2.3} Since $T$ is affine, the vector
bundle $\D_T^k$ is globally generated. Choose finitely many
differential operators
$$D_{\alpha}\in \Gamma (T, \D^k_T) $$
which span $\D^k_T$ at every point of $T$. Let $\sigma \in \Gamma(X
\times T, pr_1^*(L))$ be the section defining the given divisor $E$,
and consider the algebraic subset
$$X \times T \supset B = \big \{ (x,t) \bigm | D_\alpha \sigma (x, t)
= 0 \ \ \forall \ \alpha \big \}$$
 cut out by the common zeroes of all the sections
$D_{\alpha}\sigma \in \Gamma(X \times T, pr_1^*(L))$.

We assert that
$$V \nsubseteq B . \tag *$$
To verify this, we study the first projection
$$pr_1 : X \times T \lra X.$$
Fix any point $x \in X$, and consider the fibre $E_x \subset T$ of $E$
over $x$. Assume that $E_x \ne T$ (which will certainly hold for
general $x$), so that $E_x$ is a divisor on $T$. Given $t \in T$, it
follows from the fact that the $D_{\alpha}$ generate $\D^k_T$ at $t$
that
$$(x, t) \in B \ \ \iff \ \ \mult_t(E_x) > k.$$
On the other hand, since $V$ dominates $X$, Lemma 2.1 applies to $pr_1$
and we conclude that
$$\mult_t(E_x) = \mult_V(E) = k$$
for sufficiently general $(x, t) \in V$. This proves (*).

It follows from (*) that if $D \in \Gamma (  T, \D^k_T)$ is a
sufficiently general $\C$-linear combination of the $D_{\alpha}$, then
$$ \sigma^\pr =_{\text{def}} D \sigma \in \Gamma( X \times T,
pr^*_1(L))$$
does not vanish on $V$. On the other hand, a differential operator of
order $\le k$ decreases multiplicities by at most $k$. Therefore if
$E^\pr$ is the divisor of $\sigma^\pr$, then $\mult_Z(E^\pr) \ge \ell -
k$, as required. \qed \enddemo

\vskip 8pt
\ni{\bf \S3. The Main Theorem.}

The purpose of  this section is to  prove:
\proclaim{Theorem 3.1}
Let $L$ be a nef  line bundle on an $n$-dimensional irreducible
projective variety $X$. Suppose there exists a
countable union $ \B \subset X$ of proper subvarieties of $X$
plus a positive real number $\a>0$ such that
$$(L)^{r}\cv Y \ge (r\cdot \a)^r \tag 3.1.1$$
for every irreducible subvariety $Y \subset X$ of dimension $r$ $(1 \le r
\le n)$ with $Y \not \subseteq \B$. Then
$$ \eps (L,x)\ge  \a$$
for all $x\in X$ outside the union of countably many proper
subvarieties of $X$.
\endproclaim
\ni Observe that (3.1.1) implies that $L$ is big. Recall also that
a line bundle $B$ is  big if and only if  there exists an ample
divisor $A$ and an effective divisor $E$ such that $aB = A + E$ for
some $a \gg 0$ (cf. \cite{Mo},(1.9)). Hence given a nef and big line
bundle $L$ on $X$, the restriction of $L$ to $Y \not \subset E$ is
again  big, and hence the inequality (3.1.1) automatically holds with
$\a = \frac{1}{n}$.  Therefore (3.1)   implies the first statement of
Theorem 1 from the Introduction, and Lemma 1.4 yields the second
assertion. Similarly, Corollary 2 follows from (3.1) and the second
statement in (1.3.1). We will prove Corollary 3 in
\S4.

(3.2). Turning to the proof of Theorem 3.1, we start with some preliminary
remarks and reductions.  First, the statement is clear if dim $X = 1$.
Therefore we may -- and do -- assume inductively that the
Theorem is known for all varieties of dimension $< n$.

Note next that there is no loss of generality in supposing that $X$ is
smooth.  In fact, let $$ f : X^\prime \lra X$$ be a resolution of
singularities,  and set $L^\prime = f^*L$, so that $L^\prime$ is a nef
line bundle  on $X^\prime$. Suppose that $Y\subset X$ is an
$r$-dimensional subvariety of
$X$, not contained in the fundamental locus of $f$. If $Y^\prime \subset
X^\prime$ is the proper transform of $Y$, then
$(L)^r\cv Y=(L^\prime)^r\cv Y^\prime$. Simililarly, if $x \in X$ is a point
over which $f$ is an  isomorphism, then one sees from (1.1) that
 $$\eps(L,x) = \eps(L^\prime, x^\prime).$$  Thus it suffices to prove the
theorem for $X^\prime$, so we will henceforth assume that $X$ is smooth.

(3.3).  Let $\beta > 0$ be a  real number and let $x \in X$ be a point at
which  $\eps(L,x) < \beta$. Then there exists a reduced irreducible curve
$C_x \subset X$ with
$$ \beta\cdot\mult_x(C_x) >  (L \cdot C_x).$$
Observe that the set of all pairs
$$ \bigg \{ (C, x) \biggm | C \subset X \text{ an integral curve},
\ \beta\cdot\mult_x(C) > (C \cdot L) \bigg \}$$
is parametrized by countably many irreducible quasi-projective
varieties. This is a consequence of the existence of Hilbert schemes,
plus the fact that in a flat family of curves, it is a constructible
condition to be reduced and irreducible (cf. \cite{Jo}, (4.10)).
It follows to begin with that the set
$$ U_{\beta}=_{\text{def}}\left \{ x \in X \biggm | \eps(L,x) < \beta \right
\}$$   can be expressed as a countable union of locally closed
subsets of $X$. Therefore to prove the Theorem, it is enough to show
that $U_{\a}$ does not contain a Zariski-open subset of $X$.
By the same token, it is even sufficient to show that for any
small rational $\d >0$  the set $U_{\a-\d}$ does not contain a
Zariski-open subset. Indeed,
$$U_{\a}=\bigcup_{\d\in \Q^+} U_{\a-\d},$$
and the latter is a countable union.

We fix  now $ \d\ll \a$ and set $\g=\a-\d$. So the issue is to show that
$U_{\g}$ does not contain a Zariski-open subset.

(3.4). Assume to the contrary that $U_\g$ does contain a Zariski-open
subset, i.e. that there exists a Zariski-open subset $U\subseteq X$ such
that
$$\eps(L, x) <\g < \a$$ for every $x \in U$. Then for every $x \in U$
there exists a reduced irreducible curve $C_x \subset X$ passing through
$x$ such that  $\g\cdot\mult_x(C_x) >  (L\cv C_x)$.  We will say that $C_x$ is
a
{\it Seshadri-exceptional} curve {\it based at} $x$.

It follows from the discussion in (3.3) that there is an irreducible
family of Seshadri exceptional curves  whose base-points sweep out an open
subset of $X$.  More precisely, there exists an
irreducible quasi-projective variety $T$, a dominant morphism
$$g : T \lra X,$$ plus an irreducible subvariety
$$ C \subset X \times T,$$ flat over $T$,
such that for every $t \in T$ the fibre
$C_t$ is a Seshadri-exceptional curve based at
$g(t) \in X$. In other words, $C_t
\subset X$ is a reduced irreducible curve, passing through $g(t)$, with
$$\g\cdot\mult_{g(t)} C_t > (L \cdot C_t).$$ Replacing $T$ first by a
suitable subvariety, and then by an open subset we can -- and do --
assume that $T$ is smooth and affine, and that $g : T \lra X$ is
quasi-finite.  Write $$\Gamma \subset X \times T$$  for the graph of $g$,
and as in (0.4) given a subset $Z \subset X\times T$, denote by $Z_t \subset
X$ the fibre of $Z$ over $t \in T$, viewed as a subset of $X$.

(3.5). We next consider a construction analogous to one used by
Koll\'ar, Miyaoka and Mori in  their proof \cite{KoMM} of the
boundedness of Fano varieties of Picard number one.
\proclaim{Lemma 3.5.1} Let
$Z \subset X \times T$ be an irreducible closed subvariety dominating both
$X$ and $T$.  Then one can construct an irreducible closed
subvariety
$$CZ \subset X \times T$$ having the following properties:
\roster
\vskip 5pt
\item"{(3.5.2)}"  $Z \subseteq CZ$ and $\text{ dim } CZ \le \text{ dim } Z +
1$.
\vskip 5pt
\item"{(3.5.3)}"  For generic $t \in T$, the fibre $(CZ)_t \subset X$
has the form
$$ (CZ)_t = \text{\rm  closure} \big (  \ \underset {s \in S_t} \to
\bigcup C_s \ \big ) ,$$
where $S_t \subset g^{-1}(Z_t)$ is a closed subset of $T$ which dominates
$Z_t$.
\endroster
\endproclaim
\ni In other words,  for general $t \in T$, $(CZ)_t$ is the closure of all
the points on a family of Seshadri exceptional curves
$\{ C_s \}_{s \in S_t}$ based at a dense constructible subset of $Z_t$.

\demo{\bf Proof}  First, let
$$S^\pr = (g \times id_T)^{-1}(Z) \subseteq T \times T.$$
The hypothesis that $Z$ dominates $X$ implies that $S^\pr \ne \emptyset$.
Fix an irreducible component $S_1$ of $S^\pr$ whose image
under $g\times id_T$ dominates $Z$. Then $\dim S_1 = \dim Z$ since $g$
is quasi-finite.  Next, letting $\pi : C \lra T$ denote the projection of
$C \subset X \times T$ onto the second factor, put
$$V_1 = (\pi \times id_T)^{-1}(S_1) \subseteq X \times T \times T.$$
Very concretely, $V_1$ may be described as the set
$$ V_1 = \big \{ (x , s, t) \bigm | x \in C_s,
\ g(s) \in Z_t, \ (s,t) \in S_1 \big \}.  $$
The fibres of the projection $p: V_1 \lra S_1$ are irreducible curves,
and hence $V_1$ is irreducible, with
$$ \dim V_1 = \dim S_1 + 1.$$
Note also that $p$ admits a section $\sigma : S_1 \lra V_1$
given by $\sigma(s,t) = (g(s), s, t)$.

Consider now the projection $pr_{13} : X \times T \times T
\lra X \times T$ onto the first and third factors, and set
$$V = pr_{13}(V_1) \subseteq X \times T.$$
Then $V$ is an irreducible constructible subset of $X \times T$,
and $V$ contains an open subset of $Z$ [viz. an open subset of
$(pr_{13} \circ \sigma)(S_1)$]. Given $t \in T$, let
$$S_t = g^{-1}(Z_t) \cap (S_1)_t,$$
where by $(S_1)_t$ we mean the fibre of $S_1\subset T\times T$ over  the
second  factor. Then by construction, for every $t \in T$:  $V_t =
\cup_{s \in S_t} C_s$.  Finally, put
$$CZ = \text{ \rm closure}(V) \subseteq X \times T.$$
Then property (3.5.2) is clear.
As for (3.5.3), it follows from the remark that if $V \subseteq X \times T$
is an irreducible constructible subset dominating $T$, then for
general $t \in T$:
$$ \text{ \rm closure}(V_t) = (\text{\rm closure}(V))_t.  $$
This completes the proof of (3.5.1). \qed \enddemo

(3.6).  The inductive hypothesis is now used to prove:
\proclaim{Lemma 3.6.1}  Let $Z \subset X \times T$ be a proper irreducible
subvariety dominating both $X$ and $T$, and consider
the variety $CZ \subseteq X
\times T$ constructed in (3.5). Then $Z$ is a proper subvariety of $CZ$.
\endproclaim
\demo{\bf Proof}  Assume to the contrary that $CZ = Z$, and fix a
very general  point $t \in T$. Given a general point $x \in Z_t$, it follows
from (3.5.3)  that there exists a Seshadri exceptional curve $C_s$ based at
$x$ such that
$C_s$ lies in $Z_t = (CZ)_t$. But this means that the restriction
$L | Z_t$ of $L$ to $Z_t$ has small Seshadri constant at a general point,
i.e.: $$\eps( L|Z_t, x) < \g$$
for a dense open set of points $x \in Z_t$. But every component of
$Z_t$ has dimension $< n$. Therefore the induction hypothesis will give a
contradiction once we show that $L | W_t$ satisfies (3.1.1) for any
irreducible component $W_t \subset Z_t$. Since the morphism $Z\lra X$ is
dominating, for sufficiently general $t\in T$ no component  $W_t$ of $Z_t$
lies entirely in $\B$. Hence for very general $t \in T$, $\B \cap W_t$ is
a countable union of proper subvarieties of $W_t$. On the other hand, if
$Y\subseteq W_t$ is a subvariety of dimension $r$ not lying in $\B \cap
W_t$, then $$(L|{W_t})^r\cv Y=L^r\cv Y \ge (r\cdot\g)^r.$$ Hence (3.1.1)
holds for $L | W_t$, as required.
\qed \enddemo

(3.7). Much as in \cite{KoMM}, (3.5.1) will be used to construct a chain
of irreducible subvarieties $Z_i \subseteq X \times T$, as follows. Start
with
$$Z_0 = \Gamma = \text{ graph}(g) \ ,\  Z_1 = C \ \subseteq X \times T,$$
and then for $1 < i \le n-1  $ apply (3.5.1) inductively to form
$$Z_{i+1} = CZ_i \subset X \times T.$$
It follows from (3.6.1) that $Z_i \subsetneq Z_{i+1}$, and
consequently $Z_i$ has relative dimension  $i$ over $T$.
In particular, $Z_n = X \times T$. Thus we have a chain
$$\Gamma = Z_0 \subset Z_1 \subset \dots \subset Z_{n-1}
\subset Z_n = X \times T
\tag 3.7.1$$
of irreducible subvarieties of $X \times T$.

(3.8). We now come to the second construction, inspired by Nadel
\cite{Nad}, Campana \cite{Ca}  and the gap arguments used in
connection with zero estimates  (cf. \cite{Be}, \cite{Fa},
\cite{Nak}).  The idea is to choose a family of divisors $E_t \in
|kL|$ $(k \gg 0)$  having high multiplicity at $g(t) \in X$, and to study
the multiplicities of $E_t$ along the  subvarieties $(Z_i)_t$ defined in
(3.7).

We start with a pointwise description. Since $(L^n) \ge (\a n)^n > (\g n)^n$,
a standard parameter count  shows that if $k \gg 0$, then
given any point
$x \in X$ there exists a divisor
$$E_x \in |kL| \ \ \text{with} \ \ \mult_x(E_x) > k\g n.$$
In fact, by Riemann-Roch
$$h^0( \O_X(kL)) = k^n \frac{(L^n)}{n!} + o(k^n),$$
whereas it is $\frac{( k \g n)^n}{n!} + o(k^n)$ conditions to impose
multiplicity $[k \gamma  n + 1]$ at a given point. In particular, we may
apply this with $x = g(t)$ to construct a divisor $E_t$ having high
multiplicity at the base of the  Seshadri-exceptional curve $C_t$.

These remarks globalize in the following manner.
Put $b = [k\g n + 1]$, and consider the projections
$pr_1 : X \times T \lra X$, \  $pr_2 : X \times T \lra T$.
Then for $k \gg 0$ the torsion-free $\O_T$-module
$$\F = pr_{2,*}(pr_1^*(kL) \otimes \I_{\Gamma}^b)$$
has positive rank, where $\I_{\Gamma} \subset \O_{X \times T}$
denotes the ideal sheaf
of $\Gamma$. As $T$ is affine, $\F$ is globally generated.
We fix a non-sero section $\sigma \in \Gamma(T, \F)$, and since
$$\Gamma(T ,\F) = \Gamma( X \times T, pr_1^*(kL) \otimes
\I_{\Gamma}^b),$$
$\sigma$ gives rise to a divisor
$$E \in |\O_{X \times T}(pr_1^*(kL))| \ \ \text{with} \ \ \mult_{\Gamma}(E) >
{k\g n}.$$

(3.9). Consider now the multiplicities $$\mult_{Z_i}(E)$$
of $E$ along the sets $Z_i$ appearing in (3.7.1). We have
$$ \mult_{Z_0}(E) = \mult_{\Gamma}(E) > k\g n\ \ , \ \
\mult_{Z_n}(E) = \mult_{(X \times T)}(E) = 0. $$
It follows that there is at least one index $i$ $(0 \le i \le n-1)$ such that
$$
\mult_{Z_i}(E) - \mult_{Z_{i+1}}(E) >  k\g. $$
The heart of the argument is that we can now apply
Propostion 2.3 to produce a new divisor, not containing $Z_{i+1}$,
with relatively high multiplicity along $Z_i$.

Specifically, since $Z_{i+1}$ dominates $X$, Proposition 2.3
implies the existence of a divisor
$$E^\pr \in | \O_{X \times T}(\ pr_1^*(kL) \ )|$$
such that
$$\mult_{Z_i}(E^\pr) > k\g \ \ \ , \  \ \ Z_{i+1} \not \subseteq
\text{Supp}(E^\pr).$$
Fix a general point $t \in T$, and consider the divisor
$E_t^\pr \in |kL|$ on $X$.
Then $E^\pr_t$ does not contain any component of $(Z_{i+1})_t$,
whereas it follows from Lemma 2.1 that
$$\mult_{W_t}(E^\pr_t) = \mult_{Z_i}(E^\pr) > k\g$$
for any irreducible component $W_t$ of $(Z_i)_t$.

Consider finally a general point $x \in W_t$ for some irreducible component
$W_t \subset (Z_i)_t$. Then $$\mult_x(E^\pr_t) > k\g.$$
On the other hand, it follows from property (3.5.3) of the construction
(3.5.1) of $Z_{i+1}$  that there is a Seshadri exceptional curve
$C_s \subset (Z_{i+1})_t$ based at $x$ such that
 $$C_s \not \subseteq \text{Supp}  (E_t^\pr).$$
Then $C_s$ meets $E^\pr_t$ properly, and we find:
$$ \align
k (L \cdot C_s) = E_t^\pr \cdot C_s &\ge \mult_x(E_t^\pr) \cdot \mult_x(C_s) \\
  &> k\g \cdot \mult_x(C_s ).
\endalign
$$
This contradicts the fact that $C_s$ is Seshadri exceptional,
and completes the proof of  Theorem 3.1.

\vskip 9pt
\ni{\bf \S 4.  Applications.}

In this section we give some simple applications of the main Theorem.

We begin with a criterion for birationality
which, together with (3.1), implies Corollary 3 from the Introduction.
\proclaim{Lemma 4.1} Let $X$ be a smooth projective
variety of dimension $n \ge 2$ and $L$ a nef and big line bundle on
$X$. Suppose that there exists a countable union $\V \subset X$ of
proper subvarieties such that $\eps(L, x) \ge 2n$ for all $x \in X -
\V$. Then the adjoint bundle
$\O_X(K_X + L)$ is birationally very ample, i.e. the corresponding rational
mapping
$$\phi_{|K_X + L|} : X  \dashrightarrow \P$$
maps $X$ birationally onto its image. \endproclaim

\demo{\bf Proof} We start with a  general remark. Suppose  that $X$ is an
irreducible projective variety and $B$ is a line bundle on
$X$, with $H^0(X,B) \ne 0$, defining a rational mapping
$$\phi = \phi_{|B|} : X  \dashrightarrow \P.$$
Then we claim that $\phi$ is birational onto its image if and only if
there exists a countable union  $\V \subset X$ of proper subvarieties
such that $\phi$ is defined and one-to-one on $X - \V$. In fact,
there exists in any event a Zariski-open subset $U \subset X$
(which in general may be empty) on which $\phi$ is defined and
one-to-one. The stated condition implies that $U \ne \emptyset$, and
so $\phi$ is generically one-to-one over its image, hence birational.

Returning to the situation of the Lemma, take $B = K_X + L$. We
will prove momentarily that for any two distinct points $x , y \in X -
\V$ one has the the vanishing
$$H^1(X, \O_X(K_X + L) \otimes \I_x \otimes \I_y) = 0.
\tag *$$
But this means exactly that $\phi_{|K_X + L|}$ is defined and
one-to-one on $X - \V$, and hence is birational onto its image, as
claimed. As for (*), let
$f : Y \lra X$ be the blowing up of $X$ at $x$ and $y$, and denote by $E_x ,
E_y \subset Y$ the
exceptional divisors.  Since $\eps(L, x) \ge 2n$ and $\eps(L,y) \ge
2n$ by hypothesis, it is a consequence of (1.2) that the
$\Q$-divisors
$$\frac{1}{2}f^*L - nE_x , \ \frac{1}{2}f^*L - nE_y$$
are nef.  Therefore $f^*L - nE_x - nE_y$ is nef. Moreover,  since $n
\ge 2$ we have $L^n \ge (2n)^n > 2n^n$ by (1.8), and so $f^*L - nE_x -
nE_y$ is also big. Then just as in (1.3), (*) follows from vanishing on
$Y$. \qed \enddemo

In the rest of this section we outline how these results can be
generalized in the context of $\Q$-divisors. To begin with, note:
\proclaim{Remark 4.2} Suppose that $X$ is an irreducible projective
variety, and $L$ is a nef $\Q$-Cartier $\Q$-divisor on $X$ satisfying
the numerical hypotheses (3.1.1).  Then $\eps(L,x) \ge \a$ for all
smooth $x \in X$ outside the union of countably many proper
subvarieties.\endproclaim
\ni In fact, choose a positive integer $m > 0$ such that
$mL$ is a Cartier divisor. Since $\eps(mL, x) = m \cdot \eps(L,x)$ for
all $x \in X$,  the assertion follows from (3.1).

We will henceforth deal with the following set-up:
\subhead
Assumptions 4.3
\endsubhead
 $X$ is a smooth  irreducible projective variety of dimension $n\ge
2$, and $L$ is a nef and big  $\Q$-divisor on $X$.  We suppose
that  $\Delta$ is a fractional  $\Q$-divisor on $X$ (i.e.
$\llcorner\Delta\lrcorner= 0$)  with normal crossing support.
Finally we assume that $N$ is an integral divisor on $X$ satisfying
the numerical equivalence $N\equiv L+\Delta$.

Arguing much as in the proof of (1.3.1), but using Kawamata-Viehweg
vanishing for $\Q$-divisors, one then finds first of all:
\proclaim{Proposition 4.4}
In the situation of  (4.3) suppose that  $L$ satisfies
the numerical hypothesis (3.1.1) of Theorem 3.1. If
$$\a \ge n + s,$$  then $|K_X+N|$ generates $s$-jets at a very general point
$x\notin\text{Supp}(\Delta)$.
\qed
\endproclaim
\ni The case $s = 0$ is proven (in more generality, and with slightly weaker
numerical hypotheses) by Koll\'ar in \cite{Ko1}, \S3. As in
\cite{Ko2}, \S8, this implies  for example that if $X$ is a smooth
projective variety with generically large algebraic fundamental group, and
$L$ is any big line bundle on $X$, then $H^0(X, \O_X(K_X + L)) \ne 0$. It
would be interesting to know whether one can use the cases $s > 0$ of
Proposition 4.4 (or the birationality statement of Proposition 4.5 below)
to obtain further information under suitable hypotheses on $L$. We note
that it follows from \cite{Ko2}, Lemma 8.2, that if $X$ has
generically large algebraic fundamental group, and $L$ is an ample line
bundle on $X$, then given any $\a > 0$ there exists an \'etale covering $m :
X^\prime \lra X$ such that $m^*L$ satisfies the hypotheses of Theorem 3.1.
However it is not immediately clear how to pass to useful information on $X$
beyond the non-vanishing established by Koll\'ar.

Arguing as in (4.1) one finds similarly:
\proclaim{Proposition 4.5}
In the set-up of (4.3), suppose that $\eps(L,x) \ge 2n$ for a very
general  point $x\in X$.
Then $\O_X(K_X+N)$ is birationally very ample.
\qed
\endproclaim
\ni In view of Remark 4.2, this applies in particular if $L$
satisfies the numerical hypotheses of (3.1.1) with $\a \ge 2n$.

Finally, we give a simple application of (4.5) to pluricanonical maps
of minimal varieties:
\proclaim{Corollary 4.6}  Let $X$ be a minimal $n$-fold of general type
having (global) index $r$, i.e. assume that $X$ has only terminal
singularities, that $K_X$ is nef and big, and that $rK_X$ is Cartier.
Then the pluricanonical series $|mrK_X|$ is birationally very ample for
$m \ge 2n^2 + 1$.
\endproclaim
\demo{\bf Sketch of Proof}  Let $f : Y \lra X$ be a log resolution of
$X$. Since $X$ has only terminal singularities, we can write
$$K_Y + \Delta \equiv f^* K_X + P,$$
where $\Delta$ is a fractional divisor
(i.e. $\llcorner\Delta\lrcorner=\emptyset$) with normal crossing support,
and $P$ is integral, effective and $f$-exceptional. Hence
$$K_Y + \Delta + (mr - 1)f^* K_X \equiv f^*(mr K_X) + P.$$
By (3.1), $\eps(f^*\O_X(rK_X), y) \ge \frac{1}{n}$ for very general
$y \in Y$. Hence if $m > 2n^2$, then
$$\eps((mr-1)f^*K_X)\ge \frac{mr-1}{nr}\ge 2n+\frac{r-1}{nr},$$
and it follows from (4.5) applied to $N\equiv \Delta
+ (mr-1)f^*K_X$ that the linear
series
$|f^*(mrK_X) + P|$ is birationally very ample on $Y$. But since $P$ is
$f$-exceptional,
$$H^0(Y, \O_Y(f^*(mrK_X) + P)) = H^0(Y, \O_Y(f^*(mrK_X))) = H^0(X,
\O_X(mrK_X)).$$
Therefore $|mrK_X|$ is birationally very ample on $X$. \qed \enddemo

\bl
\vskip 10pt

\ni 

\vskip 20pt
\bl
\settabs\+University of Illinois at Chicago and now is the time  \cr
\+ Lawrence EIN \cr
\+ Department of Mathematics \cr
\+ University of Illinois at Chicago \cr
\+ Chicago, IL  60680 \cr
\+ e-mail: U22425$\%$UICVM.BITNET \cr

\vskip 7 pt

\+ Oliver K\"UCHLE and Robert LAZARSFELD \cr
\+ Department of Mathematics \cr
\+ University of California, Los Angeles \cr
\+ Los Angeles, CA  90024 \cr
\+ e-mail:  kuechle$\@$math.ucla.edu  \ \ and  \ \ rkl$\@$math.ucla.edu \cr
\enddocument
\bye